\makeatletter

\documentclass{paper}
\usepackage[utf8]{inputenc}
\usepackage{amsmath}
\usepackage{amsthm}
\usepackage{amssymb}
\usepackage{graphicx}
\usepackage[numbers]{natbib}
\usepackage[unicode=true,
 bookmarks=false,
 breaklinks=false,pdfborder={0 0 1},backref=section,colorlinks=false]
 {hyperref}

\makeatletter
\newenvironment{lyxcode}
	{\par\begin{list}{}{
		\setlength{\rightmargin}{\leftmargin}
		\setlength{\listparindent}{0pt}
		\raggedright
		\setlength{\itemsep}{0pt}
		\setlength{\parsep}{0pt}
		\normalfont\ttfamily}%
	 \item[]}
	{\end{list}}

\newtheorem{Remark}{Remark}

\usepackage{stmaryrd}
\makeatother

\begin{document}
\title{Modeling and control of the rodwheel}
\author{Luc Jaulin}
\institution{Lab-Sticc, ENSTA-Bretagne}

\maketitle
\textbf{Abstract}. The rodwheel is a wheel equipped with a rod motorized
on the axle. This paper proposes a Lagrangian approach to find the
state equations of the rodwheel rolling on a plane without friction.
The approach takes advantage of a symbolic computation. A controller
is proposed to stabilize the rodwheel with the rod upward and going
straight at a desired speed. 

\section{Introduction}

\begin{figure}[h]
\centering\includegraphics[width=9cm]{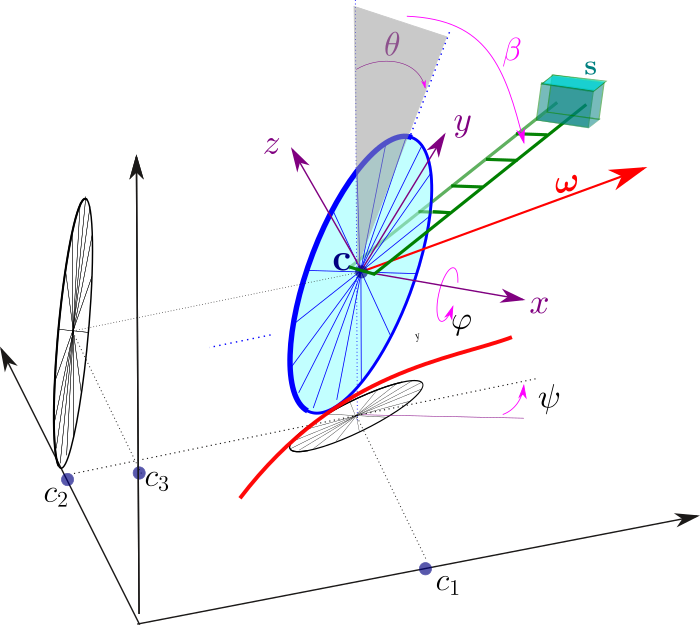}\caption{Disk rolling on a plane with a rod$ $}

\label{fig:rollingrodwheel}
\end{figure}

Consider a disk rolling on a plane without friction with a rod which
can move along the wheel plane \citep{Esnault_rapport2023} as shown
on Figure \ref{fig:rollingrodwheel}. In the axle of the wheel, we
have a motor which produces a torque $u$ between the wheel and the
rod.

The disk has a mass $m=5\text{kg}$ and its radius is $r=1m$. The
gravity is taken as $g=9.81ms^{-2}.$ The mass $\mu=1\text{kg}$ of
the rod is reduced to the single point $\textbf{s}$ which is at a
distance $\ell=2m$ to the center $\mathbf{c}$ of the wheel.

In this paper, we want to find the state equations describing the
motion of the rodwheel. The rolling wheel model has already been found
since a long time (see \emph{e.g.} \citep{Appel:1900}) and even highly
studied since (see \emph{e.g.} \citep{OReilly:disk}). Extension to
more general wheel based vehicle have also been proposed (see \emph{e.g.}
\citep{BoyerMauny18} for the bicycle). The problem of the wheel with
the rod seems to be original and I believe it can be used as a fast,
light terrestrial robot to move in a urban environment with a reduced
amount of batteries. Computing the state equations for a rodwheel
is a tedious task. In this paper, I want to take advantage of symbolic
computing (here the \texttt{sympy} package of \texttt{Python}) in
order to derive these state equations. The Lagrangian approach \citep{wells67},
often applied to model robots \citep{Corke11}, will be chosen. 

Once we have found the state model of the rodwheel, we will propose
a controller to make the rodwheel going straight ahead at a desired
speed with the rod upward. The control approach uses techniques dedicated
to underactuated mechanical systems \citep{fantoni:02}. 

\section{Modeling}

\subsection{State vector}

We take the state vector $\mathbf{x}=(c_{1},c_{2},\varphi,\theta,\psi,\beta,\dot{\varphi},\dot{\theta},\dot{\psi},\dot{\beta})$
where $(c_{1},c_{2})$ is the vertical projection of center $\mathbf{c}=(c_{1},c_{2},c_{3})$
of the disk, $\varphi,\theta,\psi$ are the three Euler angles and
$\beta$ is the rod angle. As illustrated by Figure \ref{fig:rollingrodwheel},
\begin{itemize}
\item $\varphi$ is the spin angle
\item $\theta$ is the stand angle, \emph{i.e.}, when $\theta=0$, the disk
is vertical
\item $\psi$ is the heading, \emph{i.e.}, the horizontal orientation of
the disk
\item $\beta$ is the angle of the rod with respect to the vertical axis.
\end{itemize}
Therefore, we should be able to find a state equations of the form
\begin{equation}
\dot{\mathbf{x}}=\mathbf{f}(\mathbf{x},u).\label{eq:fxu}
\end{equation}
to represent the dynamics of the rodwheel. For this, we will follow
the approach presented in \citep{jaulin:rollingdisk}.

To use \texttt{sympy}, we first declare the symbolic variables and
functions to be used:
\begin{lyxcode}
from~sympy~import~{*}

t=symbols('t')~~~~~

m,g,r,$\mu$,l=~symbols('m~g~r,~$\mu$,l')~~~~~

c1,c2~=~Function('c1')('t'),Function('c2')('t')~~~~~

dc1,dc2~=~Function('dc1')('t'),Function('dc2')('t')~~~~~

ddc1,ddc2~=~Function('ddc1')('t'),Function('ddc2')('t'),~~~~~

$\varphi$,$\theta$,$\psi$,$\beta$~=~Function('$\varphi$')('t'),Function('$\theta$')('t'),Function('$\psi$')('t'),Function('$\beta$')('t')~~~~~

d$\varphi$,d$\theta$,d$\psi$,d$\beta$~=~Function('d$\varphi$')('t'),Function('d$\theta$')('t'),

~~~~~~~~~~~~~~~Function('d$\psi$')('t'),Function('d$\beta$')('t')~~~~~

dd$\varphi$,dd$\theta$,dd$\psi$,dd$\beta$~=~Function('dd$\varphi$')('t'),Function('dd$\theta$')('t'),

~~~~~~~~~~~~~~~~~~~Function('dd$\psi$')('t'),Function('dd$\beta$')('t')

$\lambda$1,$\lambda$2~=~Function('$\lambda$1')('t'),Function('$\lambda$2')('t')
\end{lyxcode}

\subsection{Orientation}

The orientation of the disk is fixed by the three Euler angles $\varphi,\theta,\psi$.
The corresponding orientation matrix is

\begin{equation}
\mathbf{R}_{\text{euler}}(\varphi,\theta,\psi)=\left(\begin{array}{ccccc}
\cos\psi &  & -\sin\psi &  & 0\\
\sin\psi &  & \cos\psi &  & 0\\
0 &  & 0 &  & 1
\end{array}\right)\mbox{\ensuremath{\cdot}}\left(\begin{array}{c}
\cos\theta\\
0\\
-\sin\theta
\end{array}\right.\begin{array}{c}
0\\
1\\
0
\end{array}\left.\begin{array}{c}
\sin\theta\\
0\\
\cos\theta
\end{array}\right)\cdot\left(\begin{array}{ccccc}
1 &  & 0 &  & 0\\
0 &  & \cos\varphi &  & -\sin\varphi\\
0 &  & \sin\varphi &  & \cos\varphi
\end{array}\right)
\end{equation}

It is built by the following \texttt{sympy} function
\begin{lyxcode}
def~Reuler($\varphi,\theta,\psi$):

~~~R$\varphi$~=~Matrix({[}{[}1,0,0{]},{[}0,cos($\varphi$),-sin($\varphi$){]},{[}0,sin($\varphi$),cos($\varphi$){]}{]})

~~~R$\theta$~=~Matrix({[}{[}cos($\theta$),0,sin($\theta$){]},{[}0,1,0{]},{[}-sin($\theta$),0,cos($\theta$){]}{]})~~~~~~~~~

~~~R$\psi$~=~Matrix({[}{[}cos($\psi$),-sin($\psi$),0{]},{[}sin($\psi$),cos($\psi$),0{]},{[}0,0,1{]}{]})~~~~~~~~~

~~~Return~~R$\psi${*}R$\theta${*}R$\varphi$
\end{lyxcode}
The rotation vector depends on the Euler angles and their derivatives.
Its expression \citep{jaulinISTEroben} can be obtained using the
relation 
\begin{equation}
\mathbf{R}^{\text{T}}\mathbf{\dot{R}}=\left(\begin{array}{ccc}
0 & -\omega_{r3} & \omega_{r2}\\
\omega_{r3} & 0 & -\omega_{r1}\\
-\omega_{r2} & \omega_{r1} & 0
\end{array}\right)
\end{equation}
which gives us the following \texttt{sympy} function
\begin{lyxcode}
def~wr(R):~

~~~W=Transpose(R){*}diff(R,t)~

~~~return~Matrix({[}{[}-W{[}1,2{]}{]},{[}W{[}0,2{]}{]},{[}-W{[}0,1{]}{]}{]})
\end{lyxcode}

\subsection{Lagrangian}

In order, to get the state equations the rodwheel, we use a Lagrangian
approach. For this, we need to express the Lagrangian $\mathcal{L}$
with respect to the state variables. Recall that 
\begin{equation}
\mathcal{L}=E_{K}-E_{p}
\end{equation}
where $E_{K}$ is the kinetic energy and $E_{p}$ is the potential
energy. We have

\begin{equation}
E_{K}=\frac{1}{2}\boldsymbol{\omega}_{r}^{\text{T}}\mathbf{I}\boldsymbol{\omega}_{r}+\frac{1}{2}m\|\dot{\mathbf{c}}\|^{2}+\frac{1}{2}\mu\|\dot{\mathbf{s}}\|^{2}
\end{equation}
where $\mathbf{c}$ is the center of the wheel, $\mathbf{s}$ is the
end point of the rod given by 
\[
\mathbf{s}=\mathbf{c}+\mathbf{R}_{\text{euler}}(\beta,\theta,\psi)
\]
 and 
\begin{equation}
\mathbf{I}=\left(\begin{array}{lll}
\frac{mr^{2}}{2} & 0 & 0\\
0 & \frac{mr^{2}}{4} & 0\\
0 & 0 & \frac{mr^{2}}{4}
\end{array}\right)
\end{equation}
is the inertia matrix of the disk. Denote by
\begin{equation}
\mathbf{q}=\left(c_{1},c_{2},\varphi,\theta,\psi,\beta\right)
\end{equation}
the generalized coordinates of the system, \emph{i.e.}, the degrees
of freedom. The Lagrangian, which is a function of $(\mathbf{q},\dot{\mathbf{q}})$,
is computed with \texttt{sympy} as follows:
\begin{lyxcode}
def~Lagrangian(q,dq):~~~~~

~~~c1,c2,$\varphi$,$\theta$,$\psi,\beta$=list(q)~~~~~

~~~R~=~Reuler($\varphi$,$\theta$,$\psi$)~~~~~

~~~c=Matrix({[}{[}c1{]},{[}c2{]},{[}r{*}cos($\theta$){]}{]})~~~~~

~~~dc=diff(c,t)~~~~~

~~~s=Reuler($\beta,\theta,\psi$){*}Matrix({[}{[}0{]},{[}0{]},{[}l{]}{]})+c~~~~~

~~~ds=diff(s,t)~~~~~

~~~Ep=m{*}g{*}c3+$\mu${*}s{[}2{]}~~~~~

~~~I=Matrix({[}{[}1/2{*}m{*}r{*}{*}2,0,0{]},{[}0,1/4{*}m{*}r{*}{*}2,0{]},{[}0,0,1/4{*}m{*}r{*}{*}2{]}{]})~~~~~

~~~Ek=1/2{*}m{*}(dc.dot(dc))+1/2{*}$\mu${*}(ds.dot(ds))+1/2{*}wr(R).dot(I{*}wr(R))~~~~~

~~~L=subsdiff(Ek-Ep)

~~~return~L
\end{lyxcode}
where \texttt{subsdiff(E)} transforms the differentiation operator
into the corresponding state variables:
\begin{lyxcode}
def~subsdiff(E):~~~~~

~~~E=E.subs(\{diff(dc1,t):~ddc1,diff(dc2,t):~ddc2,~~diff(d$\varphi$,t):~dd$\varphi$,~

~~~~~~~~~~~~~diff(d$\theta$,t):~dd$\theta$,~diff(d$\psi$,t):~dd$\psi$,~diff(d$\beta$,t):~dd$\beta$\})

~~~E=E.subs(\{diff(c1,t):~dc1,diff(c2,t):~dc2,~~diff($\varphi$,t):~d$\varphi$,~

~~~~~~~~~~~~~diff($\theta$,t):~d$\theta$,diff($\psi$,t):~d$\psi$,~diff($\beta$,t):~d$\beta$\})~~~~~

~~~return~simplify(E)

~~~~~~~~~
\end{lyxcode}
We get
\[
\begin{array}{ccc}
\mathcal{L}(\mathbf{q},\dot{\mathbf{q}}) & = & \frac{m}{8}r^{2}(\dot{\theta}\sin\varphi-\dot{\psi}\cos\theta\cos\varphi)^{2}+\frac{m}{8}r^{2}(\dot{\theta}\cos\varphi+\dot{\psi}\sin\varphi\cos\theta)^{2}\\
 &  & +\frac{m}{4}r^{2}(\dot{\varphi}-\dot{\psi}\sin\theta)^{2}+\frac{m}{2}(r^{2}\dot{\theta}^{2}\sin^{2}\theta+\dot{c}_{1}^{2}+\dot{c}_{2}^{2})\\
 &  & -\mu(\ell\cos\beta+r)\cos\theta\\
 &  & +\frac{\mu}{2}\left(\ell\sin\psi\cos\beta\cos\theta\dot{\theta}+\ell\left(\sin\beta\sin\psi+\sin\theta\cos\beta\cos\psi\right)\dot{\psi}-\ell(\sin\beta\sin\theta\sin\psi+\cos\beta\cos\psi)\dot{\beta}+\dot{c}_{2}\right)^{2}\\
 &  & +\frac{\mu}{2}\left(\ell\cos\beta\cos\theta\cos\psi\dot{\theta}+\ell\left(\sin\beta\cos\psi-\sin\theta\sin\psi\cos\beta\right)\dot{\psi}+\ell\left(\sin\psi\cos\beta-\sin\beta\sin\theta\cos\psi\right)\dot{\beta}+\dot{c}_{1}\right){}^{2}\\
 &  & +\frac{\mu}{2}\left(\ell\sin\beta\cos\theta\dot{\beta}+(\ell\sin\theta\cos\beta+r\sin\theta)\dot{\theta}\right)^{2}\\
 &  & -gmr\cos\theta
\end{array}
\]
In what follows, we will avoid writing such complex expressions (which
can be far larger than this one). Indeed, these expressions should
remain transparent for the user since is it handled by \texttt{sympy}. 

\subsection{Euler-Lagrange equations}

The evolution of $\mathbf{q}$ obeys to the Euler-Lagrange equation
for non-holonomic systems:
\begin{equation}
\underset{\mathcal{Q}(\mathbf{q},\dot{\mathbf{q}},\ddot{\mathbf{q}})}{\underbrace{\frac{d}{dt}\left(\frac{\partial\mathcal{L}}{\partial\dot{\mathbf{q}}}\right)-\frac{\partial\mathcal{L}}{\partial\mathbf{q}}}}=\underset{\boldsymbol{\tau}}{\underbrace{\left(\begin{array}{c}
\tau_{\mathbf{c}_{1}}\\
\tau_{\mathbf{c}_{2}}\\
\tau_{\varphi}\\
\tau_{\theta}\\
\tau_{\psi}\\
\tau_{\beta}
\end{array}\right)}}+\left(\begin{array}{c}
0\\
0\\
1\\
0\\
0\\
-1
\end{array}\right)\cdot u\label{eq:EulerLagEq}
\end{equation}
where the right hand side corresponds the generalized forces. It is
composed of the forces $\boldsymbol{\tau}$ generated by the ground
and of the motor forces corresponding to internal torque. The entries
for the motor force vector are consistent with the fact that the motor
produces the torque $u$ between the degrees of freedom $\beta$ (the
rod) and $\varphi$ (the wheel). 
\begin{Remark}
\label{rem:Q:lin:ddq}Note that $\mathcal{Q}(\mathbf{q},\dot{\mathbf{q}},\ddot{\mathbf{q}})$
is affine in $\ddot{\mathbf{q}}$. This property will be used later
to build the state space model of the rodwheel. 
\end{Remark}

An expression for $\mathcal{Q}(\mathbf{q},\dot{\mathbf{q}},\ddot{\mathbf{q}})$
is obtained using \texttt{sympy} by: 

\begin{lyxcode}
q=Matrix({[}c1,c2,$\varphi$,$\theta$,$\psi$,$\beta${]})~~~~~

dq=Matrix({[}dc1,dc2,d$\varphi$,d$\theta$,d$\psi$,d$\beta${]})~~~~~

ddq=Matrix({[}ddc1,ddc2,dd$\varphi$,dd$\theta$,dd$\psi$,dd$\beta${]})

L=Lagrangian(q,dq)~

Q=subsdiff(diff(L.jacobian(dq),t)-L.jacobian(q))
\end{lyxcode}
The generalized forces $\tau_{\mathbf{c}_{1}},\tau_{\mathbf{c}_{2}},\tau_{\varphi},\tau_{\theta},\tau_{\psi},\tau_{\beta}$
generated by the reaction of the ground are linked by constraints
that should now be written.

\subsection{Non holonomic constraints}

\begin{figure}[h]
\centering\includegraphics[width=12cm]{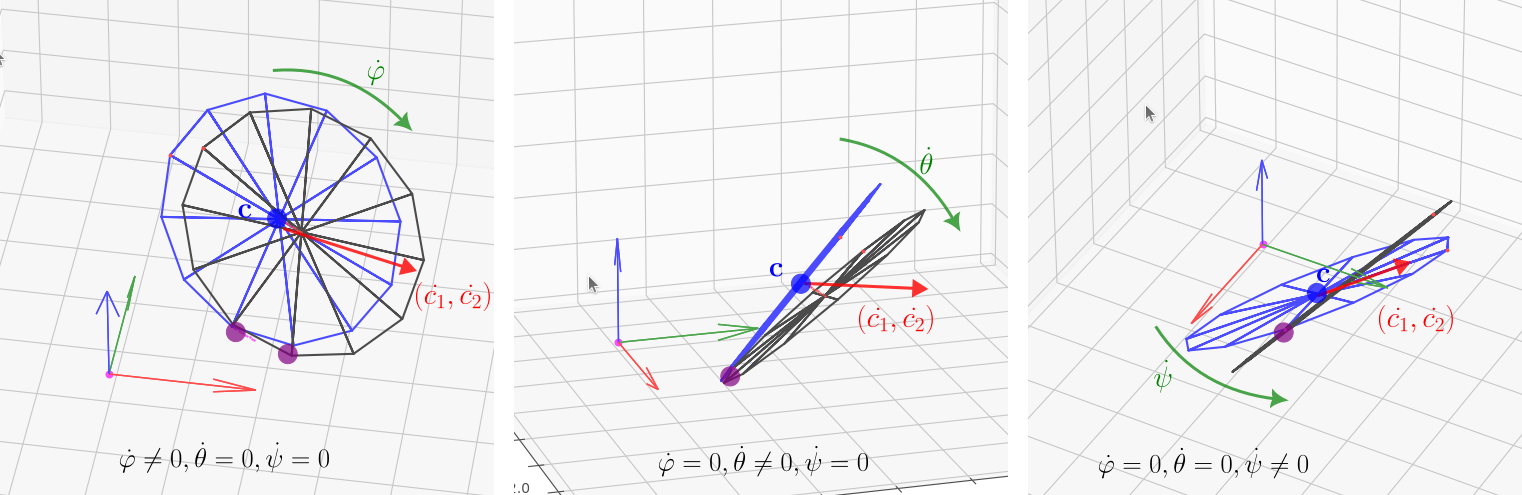}\caption{Wheel rolling on a plane (the rod is not represented)}

\label{fig:dalembert.png}
\end{figure}

Without any reaction constraint of the ground onto the disk, the state
vector would have been 
\begin{equation}
(\mathbf{q},\dot{\mathbf{q}})=(c_{1},c_{2},\varphi,\theta,\psi,\beta,\dot{c}_{1},\dot{c}_{2},\dot{\varphi},\dot{\theta},\dot{\psi},\dot{\beta}),
\end{equation}
\emph{i.e.}, it would have been composed of the degrees of freedom
$\mathbf{q}$ and they derivatives $\dot{\mathbf{q}}$. Now, due to
the ground forces, the variables $\text{(\ensuremath{\mathbf{q}},\ensuremath{\dot{\mathbf{q}}})}$
are linked by some differential constraints. These constraints will
be needed to derive the state equations with $\mathbf{x}=(c_{1},c_{2},\varphi,\theta,\psi,\beta,\dot{\varphi},\dot{\theta},\dot{\psi},\dot{\beta})$
as a state vector. Since we have two variables to eliminate (here
$\dot{c}_{1}$ and $\dot{c}_{2}$), we need to find two more differential
constraints. These two constraints are those generated by the ground
forces. They translate the fact that the point of the disk in contact
with the ground has a zero velocity. It means that the disk neither
slides tangentially (it could be the first equation) nor laterally
(it could be the second equation). We understand that these two equations
have the form
\begin{equation}
\begin{array}{ccc}
\dot{c}_{1} & = & \alpha_{1}\cdot\dot{\varphi}+\alpha_{2}\cdot\dot{\theta}+\alpha_{3}\cdot\dot{\psi}\\
\dot{c}_{2} & = & \beta_{1}\cdot\dot{\varphi}+\beta_{2}\cdot\dot{\theta}+\beta_{3}\dot{\cdot\psi}
\end{array}
\end{equation}
 where the $\alpha_{i}$'s and the $\beta_{i}'s$ depend on $\mathbf{q}$.
More precisely, as shown in \citep{jaulin:rollingdisk} these constraints
are
\begin{equation}
\begin{array}{ccc}
\dot{c}_{1} & = & r\sin\psi\cdot\dot{\varphi}+r\cos\psi\cos\theta\cdot\dot{\theta}-r\sin\psi\sin\theta\cdot\dot{\psi}\\
\dot{c}_{2} & = & -r\cos\psi\cdot\dot{\varphi}+r\sin\psi\cos\theta\cdot\dot{\theta}+r\cos\psi\sin\theta\cdot\dot{\psi}
\end{array}\label{eq:nonholonomous:constraints}
\end{equation}
Figure \ref{fig:dalembert.png} illustrates how these formulas are
obtained . The left subfigure shows that when $\dot{\theta}=0$ and
$\dot{\psi}=0$, we have 
\begin{equation}
\begin{array}{ccc}
\dot{c}_{1} & = & r\sin\psi\cdot\dot{\varphi}\\
\dot{c}_{2} & = & -r\cos\psi\cdot\dot{\varphi}
\end{array}
\end{equation}
The subfigure in the center illustrates that if $\dot{\varphi}=0$
and $\dot{\psi}=0$, 
\begin{equation}
\begin{array}{ccc}
\dot{c}_{1} & = & r\cos\psi\cos\theta\cdot\dot{\theta}\\
\dot{c}_{2} & = & r\sin\psi\cos\theta\cdot\dot{\theta}
\end{array}
\end{equation}

The right subfigure illustrates that if $\dot{\theta}=0$ and $\dot{\varphi}=0$,
\begin{equation}
\begin{array}{ccc}
\dot{c}_{1} & = & -r\sin\psi\sin\theta\cdot\dot{\psi}\\
\dot{c}_{2} & = & r\cos\psi\sin\theta\cdot\dot{\psi}
\end{array}
\end{equation}
By superposition, we get Equation \ref{eq:nonholonomous:constraints}.
These constraints are said to be non holonomic since they will not
allow us to express our system with a state composed of some degrees
of freedom $\mathbf{q}$ and their derivatives $\mathbf{\dot{q}}$. 

\subsection{D'Alembert's principle}

We need to find an expression for $\boldsymbol{\tau}$ which occurs
in the right hand side of (\ref{eq:EulerLagEq}). The components for
$\boldsymbol{\tau}$ correspond to the generalized ground forces applied
to our system. Now, since the contact point has no velocity, the ground
force cannot modify the energy of the system. 

In order to use this information, let us to recall the principle of
d'Alembert: \emph{for arbitrary virtual displacements, the constraint
forces don't do any work}. 

The virtual displacements are infinitesimal changes $\delta\mathbf{q}=\mathbf{q}(t+dt)-\mathbf{q}(t)$
that should be consistent with some feasible trajectories $\mathbf{q}(t)$.
For our rolling system, the virtual displacements satisfy
\begin{equation}
\begin{array}{ccccc}
\delta c_{1}-r\sin\psi\cdot\delta\varphi-r\cos\psi\cos\theta\cdot\delta\theta+r\sin\psi\sin\theta\cdot\delta\psi & = & 0 & \, & (i)\\
\delta c_{2}+r\cos\psi\cdot\delta\varphi-r\sin\psi\cos\theta\cdot\delta\theta-r\cos\psi\sin\theta\cdot\delta\psi & = & 0 &  & (ii)
\end{array}\label{eq:virtual:disp}
\end{equation}
for the same reasons than those used to derive (\ref{eq:nonholonomous:constraints}).
The fact that there is no work translates into
\begin{equation}
\begin{array}{ccc}
\delta W=\boldsymbol{\tau}^{\text{T}}\cdot\delta\mathbf{q}=\tau_{c_{1}}\cdot\delta c_{1}+\tau_{c_{2}}\cdot\delta c_{2}+\tau_{\varphi}\cdot\delta\varphi+\tau_{\theta}\cdot\delta\theta+\tau_{\psi}\cdot\delta\psi & +\tau_{\beta}\cdot\delta\beta= & 0.\end{array}\label{eq:no:work}
\end{equation}

Equivalently, Equation \ref{eq:no:work} is a linear combination of
the two equations (\ref{eq:virtual:disp}), \emph{i.e.}, (\ref{eq:no:work})
$=\lambda_{1}\cdot$(\ref{eq:virtual:disp},i)$+\lambda_{2}\cdot$(\ref{eq:virtual:disp},ii).
The $\lambda_{i}$'s are called the \emph{Lagrange parameters}. Therefore:
\begin{equation}
\left(\begin{array}{c}
\tau_{c_{1}}\\
\tau_{c_{2}}\\
\tau_{\varphi}\\
\tau_{\theta}\\
\tau_{\psi}\\
\tau_{\beta}
\end{array}\right)=\lambda_{1}\cdot\left(\begin{array}{c}
1\\
0\\
-r\sin\psi\\
-r\cos\psi\cos\theta\\
r\sin\psi\sin\theta\\
0
\end{array}\right)+\lambda_{2}\cdot\left(\begin{array}{c}
0\\
1\\
r\cos\psi\\
-r\sin\psi\cos\theta\\
-r\cos\psi\sin\theta\\
0
\end{array}\right)
\end{equation}

or in a vector form:
\begin{equation}
\boldsymbol{\tau}=\mathbf{A}^{\text{T}}(\mathbf{q})\left(\begin{array}{c}
\lambda_{1}\\
\lambda_{2}
\end{array}\right)\label{eq:tau:A:lam}
\end{equation}
where
\begin{equation}
\mathbf{A}(\mathbf{q})=\left(\begin{array}{cccccc}
\,1\, & 0 & \,-r\sin\psi\, & \,-r\cos\psi\cos\theta\, & \,r\sin\psi\sin\theta\, & 0\\
0 & 1 & r\cos\psi & -r\sin\psi\cos\theta & -r\cos\psi\sin\theta & 0
\end{array}\right).
\end{equation}

\subsection{State equations}

Using (\ref{eq:EulerLagEq}) with (\ref{eq:tau:A:lam}), we get

\begin{equation}
\mathcal{Q}(\mathbf{q},\dot{\mathbf{q}},\ddot{\mathbf{q}})-\boldsymbol{\tau}(\mathbf{q},\text{\ensuremath{\boldsymbol{\lambda}}})=\left(\begin{array}{c}
0\\
0\\
u\\
0\\
0\\
-u
\end{array}\right)\label{eq:eulerlag:eq}
\end{equation}

This system is made of 6 equations which are linear in 8 variables
: $\boldsymbol{\lambda}=(\lambda_{1},\lambda_{2})$ and $\ddot{\mathbf{q}}=(\ddot{c}_{1},\ddot{c}_{2},\ddot{\varphi},\ddot{\theta},\ddot{\psi},\ddot{\beta})$
(see Equation \ref{eq:tau:A:lam} and Remark \ref{rem:Q:lin:ddq}).
In order to square the system, we need to add two equations (to get
8 equations). They are obtained by differentiating of the non-holonomic
constraints (\ref{eq:nonholonomous:constraints}) given by

\begin{equation}
\underset{\mathbf{a}(\mathbf{q},\dot{\mathbf{q}})}{\underbrace{\mathbf{A}(\mathbf{q})\cdot\dot{\mathbf{q}}}}=\mathbf{0}.
\end{equation}
 Let us differentiate this equation. We get:
\begin{equation}
\underset{=\frac{d}{dt}\mathbf{a}(\mathbf{q},\dot{\mathbf{q}})}{\underbrace{\frac{\partial\mathbf{a}(\mathbf{q},\dot{\mathbf{q}})}{\partial\mathbf{q}}\cdot\dot{\mathbf{q}}+\frac{\partial\mathbf{a}(\mathbf{q},\dot{\mathbf{q}})}{\partial\dot{\mathbf{q}}}\cdot\ddot{\mathbf{q}}}}=\mathbf{0}
\end{equation}

Adding these two equations to (\ref{eq:eulerlag:eq}) yields
\[
\left\{ \begin{array}{ccc}
\frac{d}{dt}\mathbf{a}(\mathbf{q},\dot{\mathbf{q}}) & = & \left(\begin{array}{c}
0\\
0
\end{array}\right)\\
\mathcal{Q}(\mathbf{q},\dot{\mathbf{q}},\ddot{\mathbf{q}})-\boldsymbol{\tau}(\mathbf{q},\text{\ensuremath{\boldsymbol{\lambda}}}) & = & \left(\begin{array}{c}
0\\
0\\
u\\
0\\
0\\
-u
\end{array}\right)
\end{array}\right.
\]
or equivalently
\begin{equation}
\underset{\mathcal{S}(\mathbf{q},\dot{\mathbf{q}},\ddot{\mathbf{q}},\text{\ensuremath{\boldsymbol{\lambda}}})}{\underbrace{\left(\begin{array}{c}
\frac{\partial\mathbf{a}(\mathbf{q},\dot{\mathbf{q}})}{\partial\mathbf{q}}\cdot\dot{\mathbf{q}}+\frac{\partial\mathbf{a}(\mathbf{q},\dot{\mathbf{q}})}{\partial\dot{\mathbf{q}}}\cdot\ddot{\mathbf{q}}\\
\mathcal{Q}(\mathbf{q},\dot{\mathbf{q}},\ddot{\mathbf{q}})-\boldsymbol{\tau}(\mathbf{q},\text{\ensuremath{\boldsymbol{\lambda}}})
\end{array}\right)}}=\underset{\mathbf{b}_{u}}{\underbrace{\left(\begin{array}{c}
0\\
0\\
0\\
0\\
1\\
0\\
0\\
-1
\end{array}\right)}\cdot}u\label{eq:S:lag}
\end{equation}
Since $\mathcal{S}(\mathbf{q},\dot{\mathbf{q}},\ddot{\mathbf{q}},\text{\ensuremath{\boldsymbol{\lambda}}})$
is affine in $\ensuremath{\boldsymbol{\lambda}},\ddot{\mathbf{q}}$
(see Equation \ref{eq:tau:A:lam} and Remark \ref{rem:Q:lin:ddq}),
we have
\begin{equation}
\mathcal{S}(\mathbf{q},\dot{\mathbf{q}},\ddot{\mathbf{q}},\text{\ensuremath{\boldsymbol{\lambda}}})=\underset{=\mathbf{M}(\mathbf{q})}{\underbrace{\left(\begin{array}{ccc}
\frac{\partial}{\partial\boldsymbol{\lambda}}\mathcal{S}(\mathbf{q},\dot{\mathbf{q}},\mathbf{0},\mathbf{0}) & \, & \frac{\partial}{\partial\ddot{\mathbf{q}}}\mathcal{S}(\mathbf{q},\dot{\mathbf{q}},\mathbf{0},\mathbf{0})\end{array}\right)}}\cdot\left(\begin{array}{c}
\boldsymbol{\lambda}\\
\ddot{\mathbf{q}}
\end{array}\right)+\underset{=-\mathbf{b}(\mathbf{q},\dot{\mathbf{q}})}{\underbrace{\mathcal{S}(\mathbf{q},\dot{\mathbf{q}},\mathbf{0},\mathbf{0})}}.
\end{equation}
where $\mathbf{M}(\mathbf{q})$ is a square matrix called the \emph{mass
matrix} (we can show that it does not depend on $\dot{\mathbf{q}}$,
this is why we have written $\mathbf{M}(\mathbf{q})$ instead of $\mathbf{M}(\mathbf{q},\dot{\mathbf{q}})$).
Thus (\ref{eq:S:lag}) becomes
\begin{equation}
\mathbf{M}(\mathbf{q})\cdot\left(\begin{array}{c}
\boldsymbol{\lambda}\\
\ddot{\mathbf{q}}
\end{array}\right)=\mathbf{b}(\mathbf{q},\dot{\mathbf{q}})+\mathbf{b}_{u}\cdot u.
\end{equation}

Isolating $\ddot{\mathbf{q}}$, we get 
\begin{equation}
\left(\begin{array}{c}
\ddot{\varphi}\\
\ddot{\theta}\\
\ddot{\psi}\\
\ddot{\beta}
\end{array}\right)=\underset{\mathbf{P}}{\underbrace{\left(\begin{array}{cccccccc}
0 & 0 & 0 & 0 & 1 & 0 & 0 & 0\\
0 & 0 & 0 & 0 & 0 & 1 & 0 & 0\\
0 & 0 & 0 & 0 & 0 & 0 & 1 & 0\\
0 & 0 & 0 & 0 & 0 & 0 & 0 & 1
\end{array}\right)}}\cdot\mathbf{M}^{-1}(\mathbf{q})\cdot\left(\mathbf{b}(\mathbf{q},\dot{\mathbf{q}})+\mathbf{b}_{u}\cdot u\right)\label{eq:lam:ddq}
\end{equation}
where $\mathbf{P}$ is the projection matrix which selects the four
last components of the vector $(\lambda_{1},\lambda_{2},\ddot{c}_{1},\ddot{c}_{2},\ddot{\varphi},\ddot{\theta},\ddot{\psi},\ddot{\beta}).$

The expression for $\mathbf{M}$ and $\mathbf{b}$ are obtained by:
\begin{lyxcode}
A=Matrix({[}{[}1,0,-r{*}sin($\psi$),-r{*}cos($\psi$){*}cos($\theta$),r{*}sin($\psi$){*}sin($\theta$),0{]},

~~~~~~~~~~{[}0,1,~r{*}cos($\psi$),-r{*}sin($\psi$){*}cos($\theta$),-r{*}cos($\psi$){*}sin($\theta$),0{]}{]})~~~~~

$\tau$=$\lambda$1{*}A{[}0,:{]}+$\lambda$2{*}A{[}1,:{]}~

a=A{*}dq~~~~~

da=diff(a,t)~~~~~

S=Matrix({[}da,{*}list(Q-$\tau$){]})

S=subsdiff(S)~~~~

M=S.jacobian({[}$\lambda$1,$\lambda$2,ddq{]})~~~~~

b=-S.subs(\{$\lambda$1:0,$\lambda$2:0,ddc1:0,ddc2:0,dd$\varphi$:0,~dd$\theta$:0,~dd$\psi$:0,~dd$\beta$:0\})
\end{lyxcode}
Therefore, the state equations of the rodwheel are
\begin{equation}
\begin{array}{ccccc}
\left(\begin{array}{c}
\dot{c}_{1}\\
\dot{c}_{2}
\end{array}\right) & = & r\left(\begin{array}{ccccc}
\sin\psi & \, & \cos\psi\cos\theta & \, & -\sin\psi\sin\theta\\
-\cos\psi &  & \sin\psi\cos\theta &  & \cos\psi\sin\theta
\end{array}\right)\left(\begin{array}{c}
\dot{\varphi}\\
\dot{\theta}\\
\dot{\psi}
\end{array}\right) & \,\, & \text{(see (\ref{eq:nonholonomous:constraints}))}\\
\left(\begin{array}{c}
\dot{\varphi}\\
\dot{\theta}\\
\dot{\psi}\\
\dot{\beta}
\end{array}\right) & = & \left(\begin{array}{c}
\dot{\varphi}\\
\dot{\theta}\\
\dot{\psi}\\
\dot{\beta}
\end{array}\right)\\
\left(\begin{array}{c}
\ddot{\varphi}\\
\ddot{\theta}\\
\ddot{\psi}\\
\ddot{\beta}
\end{array}\right) & = & \mathbf{P}\cdot\mathbf{M}^{-1}(\mathbf{q})\cdot\left(\mathbf{b}(\mathbf{q},\dot{\mathbf{q}})+\mathbf{b}_{u}\cdot u\right) &  & \text{(see (\ref{eq:lam:ddq}))}
\end{array}\label{eq:disk:state:eq}
\end{equation}

The symbolic computation is only used to generate the \texttt{Python}
functions for $\mathbf{M}(\mathbf{q})$ and $\mathbf{b}(\mathbf{q},\dot{\mathbf{q}})$.
Getting an expression for $\mathbf{v}(\mathbf{q},\dot{\mathbf{q}},u)=\mathbf{M}^{-1}(\mathbf{q})\cdot\left(\mathbf{b}(\mathbf{q},\dot{\mathbf{q}})+\mathbf{b}_{u}\cdot u\right)$
would have been too heavy for \texttt{sympy}. This is why the choice
has been made to compute $\mathbf{v}(\mathbf{q},\dot{\mathbf{q}},u)$
numerically during the simulation by solving the following linear
system
\[
\mathbf{M}(\mathbf{q})\cdot\mathbf{v}=\mathbf{b}(\mathbf{q},\dot{\mathbf{q}})+\mathbf{b}_{u}\cdot u
\]
each time we want to evaluate the evolution function $\mathbf{f}(\mathbf{x},u).$
To build the \texttt{Python}  function from a \texttt{sympy} expression,
we write
\begin{lyxcode}
M,b=lambdify((c1,c2,$\varphi$,$\theta$,$\psi$,$\beta$,d$\varphi$,d$\theta$,d$\psi$,d$\beta$,m,g,r,$\mu$,l),(M,b))
\end{lyxcode}

\section{Control}

We want to build a controller for the rodwheel and validate on a simulation.
We first have to build a simulator based on the state space model.

\subsection{Simulation\label{subsec:simulation}}

\begin{figure}[h]
\centering\includegraphics[width=9cm]{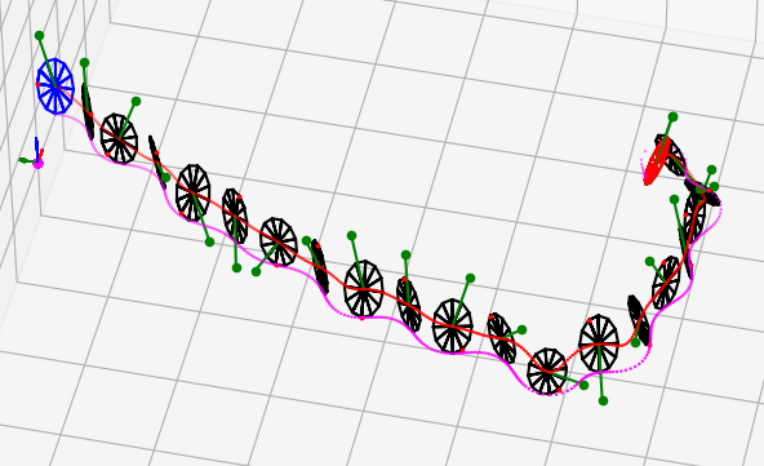}\caption{Evolution of rodwheel with no control}

\label{fig:rodweelsimu0}
\end{figure}

The \texttt{Python}  code for the evolution function is:
\begin{lyxcode}
def~f(x,u):~~~~~

~~~c1,c2,$\varphi$,$\theta$,$\psi$,$\beta$,d$\varphi$,d$\theta$,d$\psi$,d$\beta$=list(x.flatten())~~~~~

~~~dc1=~r{*}sin($\psi$){*}d$\varphi$+r{*}cos($\psi$){*}cos($\theta$){*}d$\theta$-r{*}sin($\theta$){*}sin($\psi$){*}d$\psi$~~~~~

~~~dc2=-r{*}cos($\psi$){*}d$\varphi$+r{*}sin($\psi$){*}cos($\theta$){*}d$\theta$+r{*}sin($\theta$){*}cos($\psi$){*}d$\psi$~~~~~

~~~bu=array({[}{[}0{]},{[}0{]},{[}0{]},{[}0{]},{[}1{]},{[}0{]},{[}0{]},{[}-1{]}{]})~~~~~

~~~v=np.linalg.solve(M,b+bu{*}u)~~~~

~~~dd$\varphi$,dd$\theta$,dd$\psi$,dd$\beta$=list((v{[}4:8{]}).flatten())~~~~~

~~~return~array({[}{[}dc1{]},{[}dc2{]},{[}d$\varphi${]},{[}d$\theta${]},{[}d$\psi${]},{[}d$\beta${]},{[}dd$\varphi${]},{[}dd$\theta${]},{[}dd$\psi${]},{[}dd$\beta${]}{]})
\end{lyxcode}
We use a Runge-Kutta integration scheme for the simulation 
\begin{lyxcode}
x=x+dt{*}(0.25{*}f(x,u)+0.75{*}(f(x+(2/3){*}dt{*}f(x,u),u)))~
\end{lyxcode}
We chose the initial state vector as 

\begin{equation}
\left(c_{1},c_{2},\varphi,\theta,\psi,\dot{\varphi},\dot{\theta},\dot{\psi}\right)=\left(4,0,0,0.3,0,-0.5,6,-3,0,0\right).
\end{equation}
and the parameters as 

\begin{equation}
(m,g,r,\mu,\ell)=(5,9.81,1,1,2).
\end{equation}
The result of the simulation without any control (\emph{i.e.}, for
$u=0$) is represented on Figure \ref{fig:rodweelsimu0}. The rod
is painted green. The blue wheel corresponds to $t=0$ and the red
wheel to $t=8$. The evolution of the wheel's center and of the contact
point are represented in red and magenta respectively. The simulation
highlights the chaotic behavior of the rodwheel. 

\subsection{Controller}

We want a controller that makes our rodwheel going straight, the rod
upward at a given speed. Several difficulties have to be raised:
\begin{itemize}
\item We have a single input $u$ for several state variables to control.
\item The equations of the system are complex. It is difficult to write
a closed form for the state equations in less than one page. 
\item The system is not controllable once the rodwheel goes straight.
\item The system is unstable with respect to the rod angle $\beta$.
\item The system is unstable with respect to the stand angle $\theta$.
\end{itemize}
Concerning the proof of the non controllability, instead of using
the model, it is more direct to use the symmetries of the problem.
Assume that $\dot{\psi}=0$, $\theta=0,\text{\ensuremath{\dot{\theta}=0}}$
(this is what we have once our goal is reached). We have an evolution
which remains symmetric with respect to the vertical plane supported
by the wheel. Thus, for any control $u$, we are trapped in a situation
where $\dot{\psi}=0$, $\theta=0,\text{\ensuremath{\dot{\theta}=0}}.$

For the instability a linearisation around $\dot{\psi}=0$, $\theta=0,\text{\ensuremath{\dot{\theta}=0}}$
could be used. Now, the instability with respect to $\beta$ is intuitive
and the instability with respect to $\theta$ is observed by simulation
(see further). Indeed, when the rod is stabilized above $\mathbf{c}$,
the precession tends to amplify.

In what follows, we propose a very simple expression for a possible
controller. This will be enough to illustrate that the robot can be
controlled with $u$ as a single input.

\begin{figure}[h]
\centering\includegraphics[width=9cm]{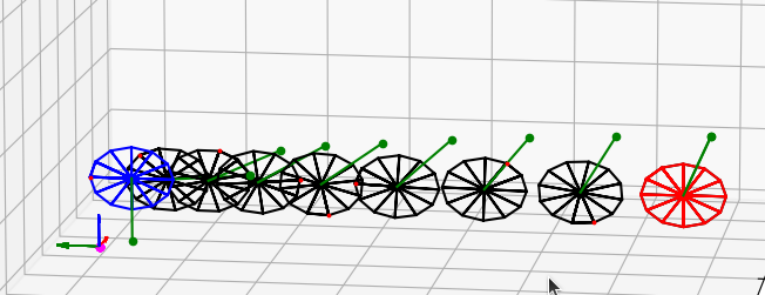}\caption{Rodwheel with $\theta(t)=0$ which accelerates to reach the speed
of $\dot{\varphi}=2$}

\label{fig:rodweelsimu1}
\end{figure}

\textbf{Case 1}. Consider the situation where the initial stand angle
$\theta$ and its derivative $\dot{\theta}$ are both equal to zero.
Due to the symmetry of the system, $\theta(t)=0$ for all $t$. Take
the following proportional and derivative controller
\begin{equation}
\begin{array}{ccc}
u & = & 20(\beta-\beta_{0})+20\dot{\beta}\\
\beta_{0} & = & \text{tanh}(2-\dot{\varphi})
\end{array}
\end{equation}
where $\beta-\beta_{0}$ is the error. If the closed loop system is
stable, then $\beta\rightarrow\beta_{0}$ (this notations means that
$\beta$ converges to $\beta_{0}$ when $t\rightarrow\infty$). Now,
the robot will accelerate (or decelerate) until $\beta_{0}=0$, \emph{i.e.},
$\dot{\varphi}=2$. This means that we can control the speed of the
rodwheel, at least when the stand angle $\theta(t)$ is equal to zero.
Figure \ref{fig:rodweelsimu1} shows the behavior of the closed loop
for the initial state :
\begin{equation}
(c_{1},c_{2},\varphi,\theta,\psi,\beta,\dot{\varphi},\dot{\theta},\dot{\psi},\dot{\beta})=(4,0,0,0,0,\pi,0,0,0,0)
\end{equation}
 For the initialization, the rod is downward and stands up in order
to create the right acceleration. We observe that $\dot{\varphi}\rightarrow2$
and that $\beta\rightarrow0$. 

\begin{figure}[h]
\centering\includegraphics[width=9cm]{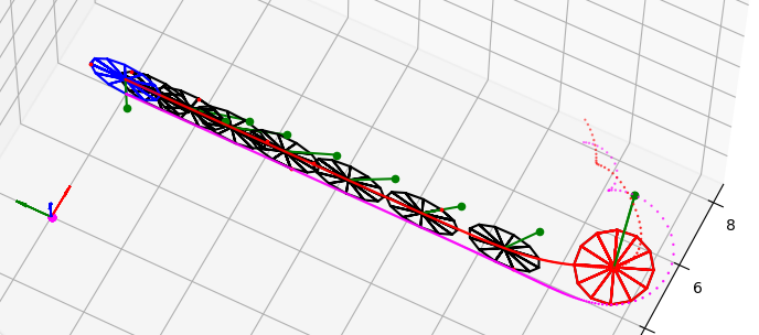}\caption{A tiny perturbation with respect to $\theta$ generates a fall of
the robot}

\label{fig:rodweelsimu1b}
\end{figure}

To show the instability with respect to $\theta$, take $\theta(0)=2\cdot10^{-12}$
instead of $\theta(0)=0$. The simulation generates Figure \ref{fig:rodweelsimu1b}.
We see the rodwheel falling strongly after few seconds.

To control the stand angle, we propose to add a term which accelerates
when $\theta$ deviates from zero. We get the following controller:
\begin{equation}
\begin{array}{ccc}
u & = & 5(\beta-\beta_{0})+5\dot{\beta}+20\cdot|\theta|\\
\beta_{0} & = & 0.2\cdot\text{tanh}(10-\dot{\varphi})
\end{array}
\end{equation}

\begin{figure}[h]
\centering\includegraphics[width=9cm]{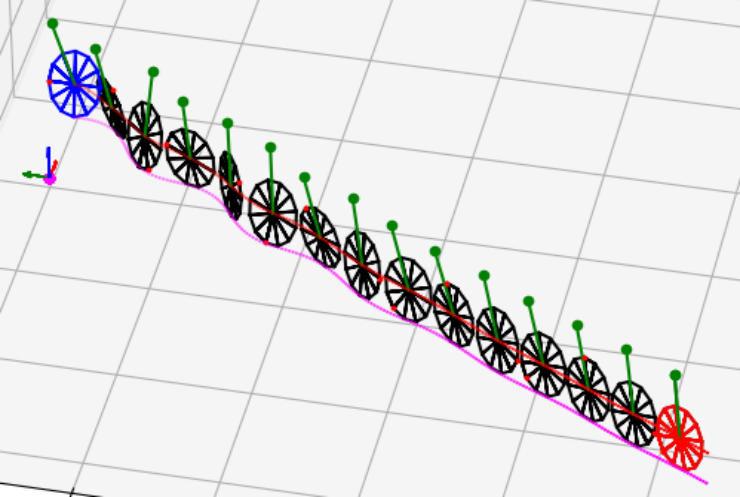}\caption{The controller stabilizes the precession of the rodwheel}

\label{fig:rodweelsimu2}
\end{figure}

The initial vector is the same as in Subsection \ref{subsec:simulation},
\emph{i.e.}, 
\begin{equation}
\left(c_{1},c_{2},\varphi,\theta,\psi,\dot{\varphi},\dot{\theta},\dot{\psi}\right)=\left(4,0,0,0.3,0,-0.5,6,-3,0,0\right).
\end{equation}
The controller is now able to limit the precession still maintaining
the rod upward as shown by Figure \ref{fig:rodweelsimu2}. 

\section{Conclusion}

In this paper, we have proposed a state space model and a simple controller
for the rodwheel. Since it is composed of one wheel, one motor and
one rod, the rodwheel can be seen as the most elementary wheeled robot
than can be built, \emph{i.e}, it seems impossible to build a fast
vehicle with less mechanical components. Due to its shape, with almost
no front surface that could slow down the vehicle, we could expect
to have a fast robot which makes long trips with few batteries. Moreover,
the expected weight could be low enough to limit the danger is case
of collision.

We have shown that it was possible to stabilize the robot and control
the speed. Now, due to the singularity for $\theta=0$ (which corresponds
to the situation where the wheel goes straight ahead) the system becomes
uncontrollable. We should thus add a other actuator such as a small
inertial wheel fixed the end of the rod to be able to control the
heading.

We have proposed a stability analysis which is based on simulation
and symmetries, which is clearly not enough. The stability should
be studied by a linearization to have a better understanding of the
instabilities and to get a better tuning of the control parameters. 

Since the rodwheel is considered here as a mobile robot, the heading
control needs to be treated. An easy solution would be to add an additional
actuated inertial wheel spinning perpendicularly with respect to the
rod. From a control point of view, a more challenging question would
be to control the heading of the robot without changing the mechanics,
\emph{i.e.} , changing the controller only. 

\bigskip{}

The \texttt{Python}  code associated to all examples can be found
here: 
\begin{center}
\href{https://www.ensta-bretagne.fr/jaulin/rodwheel.html}{https://www.ensta-bretagne.fr/jaulin/rodwheel.html}
\par\end{center}

\bibliographystyle{plain}

\end{document}